\documentclass[11pt,a4paper]{article}
\usepackage{mathrsfs}
\usepackage{amsfonts}
\usepackage{}

\usepackage[leqno]{amsmath}
\usepackage{graphicx}
\usepackage{latexsym}
\usepackage{amsmath,amsfonts,amssymb,amsthm,mathrsfs,euscript,makeidx,color}
\usepackage{enumerate}
\usepackage{bbm}

\usepackage[colorlinks,linkcolor=blue,anchorcolor=green,citecolor=red]{hyperref}
\def\sqr#1#2{{\vcenter{\vbox{\hrule height.#2pt
              \hbox{\vrule width.#2pt height#1pt \kern#1pt \vrule width.#2pt}
              \hrule height.#2pt}}}}
\def\5n{\negthinspace \negthinspace \negthinspace \negthinspace \negthinspace }
\def\4n{\negthinspace \negthinspace \negthinspace \negthinspace }
\def\3n{\negthinspace \negthinspace \negthinspace }
\def\2n{\negthinspace \negthinspace }
\def\1n{\negthinspace }

\def\dbE{\mathbb{E}}
\def\dbF{\mathbb{F}}

\def\dbL{\mathbb{L}}

\def\dbP{\mathbb{P}}

\def\dbR{\mathbb{R}}
\def\dbS{\mathbb{S}}

\def\dbX{\mathbb{X}}

\def\cA{{\cal A}}

\def\cF{{\cal F}}

\def\Bx{{\bf x}}

\def\ds{\displaystyle}

\def\ns{\noalign{\ss}}
\def\no{\noindent}

\def\ss{\smallskip}
\def\ms{\medskip}

\def\q{\quad}
\def\qq{\qquad}
\def\hb{\hbox}

\def\({\Big (}
\def\){\Big )}
\def\[{\Big[}
\def\]{\Big]}

\def\lan{\langle}
\def\ran{\rangle}

\def\rf{\eqref}
\def\a{\alpha}

\def\g{\gamma}

\def\l{\lambda}

\def\si{\sigma}

\def\f{\varphi}
\def\th{\theta}
\def\o{\omega}

\def\G{\Gamma}

\def\O{\Omega}
\def\cd{\cdot}
\def\cds{\cdots}

\def\tr{\hbox{\rm tr$\,$}}

\def\pa{\partial}

\def\les{\leqslant}
\def\ges{\geqslant}
\def\bde{\begin{definition}\label}
\def\ede{\end{definition}}
\def\be{\begin{equation}}
\def\bel{\begin{equation}\label}
\def\ee{\end{equation}}
\def\bt{\begin{theorem}\label}
\def\et{\end{theorem}}
\def\bc{\begin{corollary}\label}
\def\ec{\end{corollary}}
\def\bl{\begin{lemma}\label}
\def\el{\end{lemma}}
\def\bp{\begin{proposition}\label}
\def\ep{\end{proposition}}
\def\bas{\begin{assumption}\label}
\def\eas{\end{assumption}}
\def\br{\begin{remark}\label}
\def\er{\end{remark}}
\def\bex{\begin{example}\label}
\def\ex{\end{example}}
\def\ba{\begin{array}}
\def\ea{\end{array}}
\def\ben{\begin{enumerate}}
\def\een{\end{enumerate}}
\def\square#1{\vbox{\hrule\hbox{\vrule height#1%
     \kern#1\vrule}\hrule}}
\def\rectangle#1#2{\vbox{\hrule\hbox{\vrule height#1%
     \kern#2\vrule}\hrule}}
\def\qed{\hfill \vrule height7pt width3pt depth0pt}

\font\tenbb=msbm10 \font\sevenbb=msbm7 \font\fivebb=msbm5

\newfam\bbfam
\scriptscriptfont\bbfam=\fivebb \textfont\bbfam=\tenbb
\scriptfont\bbfam=\sevenbb

\newtheorem{theorem}{\indent Theorem}[section]
\newtheorem{definition}[theorem]{\indent Definition}
\newtheorem{proposition}[theorem]{\indent Proposition}
\newtheorem{corollary}[theorem]{\indent Corollary}
\newtheorem{lemma}[theorem]{\indent Lemma}
\newtheorem{remark}[theorem]{\indent Remark}
\newtheorem{example}[theorem]{\indent Example}

\newtheorem{assumption}[theorem]{\indent Assumption}
\makeatletter
   
   \@addtoreset{equation}{section}
\makeatother

\begin{document}

\title{\bf Non-Equivalence of Stochastic Optimal Control Problems with Open and Closed Loop Controls}

\author{
Jiongmin Yong\thanks{\noindent Department of Mathematics, University of Central Florida, Orlando, FL 32816. Email: {\tt jiongmin.yong@ucf.edu}. This work is supported in part by NSF grant DMS-1812921.} ~ and ~{Jianfeng Zhang}\thanks{\noindent
Department of Mathematics, University of Southern California, Los
Angeles, CA 90089. E-mail:{\tt jianfenz@usc.} {\tt edu}. This author is
supported in part by NSF grant DMS-1908665.
}}

\date{\today}
\maketitle

{\bf Abstract.} For an optimal control problem of an It\^o's type stochastic differential equation, the control process could be taken as open-loop or closed-loop forms. In the standard literature, provided appropriate regularity, the value functions under these two types of controls are equal and are the unique (viscosity) solution to the corresponding (path-dependent) HJB equation. In this short note, {\color{black}we provide a counterexample in the path dependent setting} showing that these value functions can be different in general. 

\ms


\ms


\section{Introduction}
\label{sect-introduction}
\setcounter{equation}{0}

Consider the following controlled It\^o's type {\color{black} path dependent} stochastic differential equation (SDE, for short) over a finite time horizon $[0, T]$:
\bel{Xpath}
X_t=x_0+\int_0^tb(s,X_{[0, s]},\a_s)ds+\int_0^t\si(s,X_{[0,s]},\a_s)dB_s,\qq t\in[0,T];
\ee
with utility functional
\bel{J}
J(\a):=\dbE\big[g(X_{[0,T]})\big].
\ee
Here $B$ is a $d$-dimensional standard Brownian motion; the controlled {\it state process} $X_t$ takes values in $\dbR^n$; $X_{[0, s]}$ refers to the path of $X$ on $[0, s]$; the coefficients $b, \si, g$ are deterministic measurable functions with appropriate dimensions, in particular $g$ is scalar valued; the {\it admissible control} $\a\in\cA$ takes values in a subset $A$ of some Euclidean space; and we shall leave the issue of existence and/or uniqueness of the state for \rf{Xpath} to later discussions.
The {\it optimal value}, or simply the value, of the control problem is defined as:
\bel{V0state}
V_0:=\sup_{\a\in\cA}J(\a),
\ee
and we call $\a^*\in \cA$  an {\it optimal control} if $J(\a^*) = V_0$.

\ms

The value $V_0$ obviously relies on the choice of the admissible control set $\cA$. Depending on the observed information in applications, among others,  the control process $\a_t$ could be taken as the so-called {\it open-loop} or {\it closed-loop} form. An open-loop control, denoted as $\a\in \cA^o$, is such that $\a$ is $\dbF^B$-progressively measurable, while a closed-loop control, denoted as $\a\in \cA^c$, is required to be $\dbF^X$-progressively measurable. Here $\dbF^B, \dbF^X$ are the natural filtration generated by $B$ and $X$, respectively. We may define the values of the control problem accordingly:
\be
\label{V0oc}
V_0^o:=\sup_{\a\in\cA^o}J(\a),\q V_0^c :=\sup_{\a\in\cA^c}J(\a).
\ee
A natural question is: do we have
\be
\label{V=}
V^o_0=V^c_0~?
\ee
We remark that, typically it is more convenient to use {\it strong formulation} for open-loop controls and {\it weak formulation} for closed-loop controls, see Remark \ref{rem-strong} below.

\ms

In the state dependent setting: for $\Bx\in C([0, T]; \dbR^n)$,
\bel{state}
b(t, \Bx, a) = b(t, \Bx_t, a),\q \si(t, \Bx, a) = \si(t, \Bx_t, a),\q g(\Bx) = g(\Bx_T),
\ee
 the standard literature provides a positive answer to \rf{V=} by using the PDE approach, see e.g. Fleming--Soner \cite{Fleming-Soner 1992} and Yong--Zhou \cite{Yong-Zhou 1999}. Consider the following HJB equation:
 \bel{HJB}
 \left.\2n\ba{c}
\ds\pa_t v(t,x)+H(t,x,\pa_x v(t,x), \pa_{xx}^2v(t,x))=0,\q v(T,x)=g(x),\ms\\
\ds\mbox{where}\q H(t,x, z, \g) := \sup_{a\in A} \big[{1\over 2}\tr\big[ \g \si \si^\top(t,x,a)\big]+z b(t,x,a)\big].
\ea\right.
\ee
Here $(t,x,z,\g)\in[0,T]\times\dbR^n\times\dbR^{1\times n}\times\dbS^n$,  with $\dbS^n$ being the set of all $n\times n$ symmetric matrices. Then, provided that the coefficients $b, \si, g$ have appropriate regularity and the HJB equation has a unique continuous viscosity solution $v$, we  have,
\bel{V=v}
V^o_0=V^c_0=v(0,x_0).
\ee
Moreover, $v(t,x)$ is the optimal value of the control problem over $[t, T]$ with initial value $X_t = x$. The main tool for this result is the {\it dynamic programming principle} (DPP for short), from which we see that the {\it dynamic value function $v(t,x)$} (more precisely we should introduce $v^o(t,x)$ and $v^c(t,x)$) of the control problem under each  type of controls is a viscosity solution of the HJB equation and hence \rf{V=v} follows from the uniqueness of the viscosity solution.

The above result remains true in the general path dependent setting \rf{Xpath}-\rf{J}. In this case, \rf{HJB} becomes a {\it path dependent HJB equation}, or more generally a {\it path dependent PDE}, with the same Hamiltonian $H$:
 \bel{HJBpath}
\pa_t v(t, \Bx)+H(t,x,\pa_\Bx v(t,\Bx), \pa_{\Bx\Bx}^2v(t,\Bx))=0;\q v(T,\Bx)=g(\Bx).
\ee
Here $\pa_t v, \pa_\Bx v, \pa^2_{\Bx\Bx} v$ are the path derivatives of Dupire \cite{Dupire-2019}. When the equation \rf{HJBpath} has a unique continuous viscosity solution, then \rf{V=v} still holds true. We refer to Zhang \cite[Part III]{Zhang 2017} for more details of the pathwise stochastic analysis and viscosity solutions of path dependent PDEs.

We emphasize that the above arguments require the dynamic value function $v(t,x)$ or $v(t, \Bx)$ to be continuous. When $v(t,x)$ is discontinuous, although there are some nice works on discontinuous viscosity solutions, see e.g. Barles--Perthame \cite{Barles-Perthame-1987}, Barron--Jensen \cite{Barron-Jensen-1990}, Bertsch--Dal Passo--Ughi \cite{BDU}, Bardi--Capuzzo--Dolcetta \cite{BC}, Chen-Su \cite{ChenSu}, and Bertsch--Smarrazzo--Terracina--Tesei \cite{BSTT}, the theory is far from complete; especially the general uniqueness issue for the discontinuous viscosity solutions to the second order equations is still open. Consequently, we
are not able to conclude \rf{V=v} or \rf{V=} from the viscosity solution approach if the value function is discontinuous.

Our main purpose of this short note is to construct a counterexample which shows that \rf{V=} can indeed fail. This implies that, besides the practical consideration in terms of the available information, mathematically it is also crucial to choose the right type of controls, especially when the value function is discontinuous. For applications of discontinuous value functions, we refer to \cite{BC} and references cited therein. We shall remark that, for stochastic differential games, even with the desired regularity, the game values can still be very sensitive to the choice of admissible controls, see e.g.  Feinstein--Rudloff--Zhang \cite{FRZ}, Possamai--Touzi--Zhang \cite{PTZ}, and Sun--Yong \cite{Sun-Yong 2014}. We also remark that, our analysis of the values does not depend on the existence of optimal controls. Another important consequence of the failure of \rf{V=} is that an (approximately) optimal control among one type of admissible controls is not necessarily (approximately) optimal anymore among the other type  of admissible controls.

Our counterexample is constructed based on the well-known example of Tsirelson \cite{Tsirelson}, which is path dependent. Note that for the state dependent case, if a second order HJB equation is uniformly non-degenerate with continuous Hamiltonian, then, even if the terminal condition is discontinuous, the value function will become continuous for $t<T$ because the diffusion term has some effect of regularization. This is not true anymore in the path dependent case, because the regularization requires some time to take effect while the discontinuity from the terminal payoff function could be present at any time in this case.

Our counterexample is constructed in \S\ref{sect-example}. In \S\ref{sect-prelim} we formulate the problems rigorously, and in \S\ref{sect-analysis} we  provide some brief discussions on the relationship between $V^o_0$ and $V^c_0$.

\section{The Problem Formulations}
\label{sect-prelim}
\setcounter{equation}{0}
We first formulate the problems rigorously in the path dependent setting. While the counterexample will be in a simpler setting, the general formulation may clarify the concepts for non-experts and will also put the counterexample in the right perspective. Denote $\dbX_n := C([0, T]; \dbR^n)$, equipped with the uniform norm: $\|\Bx\|:= \sup_{0\les t\les T} |\Bx_t|$ for all $\Bx\in \dbX_n$. Let $A\subseteq \dbR^m$ be a proper set for the possible values of admissible controls.  Consider the path dependent SDE \rf{Xpath} with coefficients $(b, \si): [0, T] \times \dbX_n \times A\to (\dbR^n, \dbR^{n\times d})$ and $g:\dbX_n\to \dbR$.  Throughout the paper, the following assumptions will always be in force:

$\bullet$ $b, \si, g$ are bounded (for simplicity) and progressively measurable in all variables;

$\bullet$ $b, \si$ are adapted in $\Bx$ in the sense that, for $\f=b,\si$,  $\f(t,\Bx, a) = \f(t, \Bx_{[0, t]}, a)$.

\no We say the system is {\it state dependent } if \rf{state} holds.

For a filtered probability space $(\O, \cF, \dbF, \dbP)$ and a generic measurable space $E$, let $\dbL^0(\dbF, \dbP; E)$ denote the set of $E$-valued processes progressively measurable with respect to the $\dbP$-augmented filtration of $\dbF$. When $\dbP$ and/or $E$ are clear, we may omit them and simply denote the set as $\dbL^0(\dbF)$. Moreover, let $\dbF^B, \dbF^X$ denote the natural filtration generated by the Brownian motion $B$ and the state process $X$, respectively.

\begin{definition} \rm
\label{defn-weak}
(i) A {\it weak solution} of the path dependent SDE \rf{Xpath} consists of a filtered probability space $(\O, \cF, \dbF, \dbP)$ and a triplet of  processes $(B, X, \a)\in \dbL^0(\dbF, \dbP;\dbR^d\times \dbR^n\times A)$ such that $B$ is a Brownian motion under $\dbP$ and \rf{Xpath} holds true $\dbP$-a.s.

(ii) A weak solution is called a {\it strong solution} if $X$ and $\a$ are $\dbF^B$-progressively measurable.
\qed
\end{definition}

In this paper we do not discuss the existence and uniqueness of weak solutions, which requires further conditions on $b, \si$. Instead, we shall always assume the following very mild assumption:

$\bullet$ for any piecewise constant control $\a_t$ valued in $A$, SDE \rf{Xpath} admits a weak solution.

\no We now introduce the optimal values under open-loop and closed-loop controls, respectively:
\bel{V0}
\left.\ba{lll}
V^o_0 := \sup \big\{\dbE^\dbP[g(X_\cd)]: \mbox{ all weak solutions of \rf{Xpath} such that $\a\in \dbL^0(\dbF^B)$}\big\};\\
\ns V^c_0 :=  \sup \big\{\dbE^\dbP[g(X_\cd)]: \mbox{ all weak solutions  of \rf{Xpath} such that $\a\in \dbL^0(\dbF^X)$}\big\}.
\ea\right.
\ee

\begin{remark}
\label{rem-strong}
{\rm (i) For open-loop controls, under the stronger conditions that $b$ and $\si$ are uniformly Lipschitz continuous in $\Bx\in \dbX_n$, one typically uses the {\it strong formulation}. That is, we fix a probability space $(\O, \cF, \dbP)$ and a Brownian motion $B$ on it. Then for any open-loop control $\a\in \dbL^0(\dbF^B)$, the SDE \rf{Xpath} admits a unique strong  solution $X\in \dbL^0(\dbF^B)$.

(ii) For closed-loop controls, it is more convenient to use {\it weak formulation}. That is, we fix the canonical space $\O:=\dbX_{d+n}$,  the canonical processes $(B, X)$, and set $\dbF := \{\cF_t\}_{0\les t\les T} := \dbF^{B, X}$, $\cF:= \cF_T$. Then for any closed-loop control $\a\in \dbL^0(\dbF^X)$, a weak solution is mainly a probability $\dbP$ on the canonical space $\dbX_{d+n}$. We remark that, for given $\a$, there might be multiple (or no) $\dbP$ corresponding to $\a$.

(iii) For closed-loop controls, since the utility $\dbE^\dbP[g(X_\cd)]$ involves only the $\dbP$-distribution of $X$, it is quite often that we consider instead the canonical space $\O := \dbX_n$ with canonical process $X$, especially when  $\si$ is non-degenerate and hence $B$ is $\dbF^X$-progressively measurable under $\dbP$. }
\qed
\end{remark}

\begin{remark}
\label{rem-dimension}
{\rm The closed-loop control case actually includes more general situations, by increasing the dimension of the state process $X$ when needed.

(i) For the case $b = b(t, B, X, a), \si  = \si(t, B, X, a), g = g(B, X)$ and/or $\a_t = \a_t(B,X)$ depend on both $B$ and $X$, we can set $\tilde X := (B, X)$ and consider the SDE in the form of \rf{Xpath}:
\bel{XB}
d \tilde X_t = \left[\ba{c} 0\\ b(t, \tilde X, \a_t)\ea\right] dt + \left[\ba{c}I_d\\ \si(t, \tilde X, \a_t)\ea\right] dB_t.
\ee
We shall remark though, in this case the coefficients $b, \si, g$ are typically discontinuous in the $B$-component of $\tilde X$, and the PDE \rf{HJB}  or PPDE \rf{HJBpath} is always degenerate. Both features {\color{red}could} contribute to the possible discontinuity of the value function.

(ii) If we allow $\a$ to be in $\dbL^0(\dbF)$ for the general $\dbF$ in Definition \ref{defn-weak}, we may still view $\a$ as a closed-loop control by considering a further enlarged state process $\hat X := (\tilde X, \G):= (B, X, \G)$:
\bel{XBG}
d \hat X_t = \left[\ba{c} 0\\ b(t, \tilde X, \a_t)\\ \a_t \ea\right] dt + \left[\ba{c}I_d\\ \si(t, \tilde X, \a_t)\\0\ea\right] dB_t,
\ee
where $\G_t=\int_0^t\a_rdr$. Note that in this case $\a$ is always in $\dbL^0(\dbF^\G)$, and hence in $\dbL^0(\dbF^{\hat X})$.
\qed}
\end{remark}

\begin{remark}
\label{rem-SMP}
{\rm In this remark we discuss some standard approaches in the literature. These approaches require appropriate regularity conditions, which we want to avoid in this paper.

(i) For both open-loop  and closed-loop controls, under appropriate regularity conditions, the dynamic value functions  $v^o(t,\Bx)$ and $v^c(t,\Bx)$ would satisfy the {\it dynamic programming principle}, which leads to the PPDE \rf{HJBpath}.  When \rf{HJBpath} has a unique  continuous viscosity solution,  we have $v^o(t, \Bx) = v^c(t,\Bx)$ and in particular $V^o_0=V^c_0$. Moreover, from the Hamiltonian, one can construct naturally an (approximate) optimal control which is closed-loop. In particular, even for the open-loop control problem in \rf{V0}, we have closed-loop (approximate) optimal controls.

(ii) Under sufficient regularity of the coefficients, any optimal open-loop control (if it exists) would satisfy the {\it stochastic maximum principle},  a Pontryagin type maximum principle, see Peng \cite{Peng 1990} or Yong--Zhou \cite{Yong-Zhou 1999}. This method is not convenient for closed-loop control though, because it involves differentiation of the closed-loop controls $\a(t,\Bx)$ with respect to $\Bx$.  Nevertheless, the optimal open-loop control $\a^*(t, B_{[0, t]})$ obtained from the stochastic maximum principle may turn out to be $\dbF^X$-progressively measurable, and in this case we also obtain the optimal closed-loop control $\tilde \a^*(t, \Bx)$ determined by: $\tilde \a^*(t, X_{[0, t]}) = \a^*(t, B_{[0, t]})$, $\dbP$-a.s.
\qed}
\end{remark}

\section{A Counterexample}
\label{sect-example}
\setcounter{equation}{0}
In this section we construct a counterexample that $V^o_0$ and $V^c_0$ are indeed not equal.  We first recall the following well-known result of Tsirelson \cite{Tsirelson}.

\begin{lemma}
\label{lem-Tsirelson} \sl
Let $t_0:=T$ and, for $k=-1,-2,\cds$,  $t_k\downarrow0$ as $k\to -\infty$. Define
\bel{mu}
\ba{c}
\ns\ds\th(x):=x-[x]\qq\forall x\in\dbR;\q\hb{$[x]$ is the greatest integer no more than $x$},\\
\ns\ds \mu(t,\Bx):=\th\Big({\Bx_{t_k}-\Bx_{t_{k-1}}\over t_k-t_{k-1}}\Big),\qq \Bx\in \dbX_1, \q t\in [t_k, t_{k+1}),~ k\les -1.
\ea\ee
Then the following SDE has no strong solution:
\bel{Xmu}
X_t=\int_0^t \mu(s, X_\cd)ds+B_t.
\ee
\end{lemma}
We note that there is a typo in the statement of \cite[Theorem]{Tsirelson}. In the definition of the coefficient $A$ there (our $\mu$ here), the domain $t\in [t_k, t_{k-1})$  should be  $t\in [t_k, t_{k+1})$ as in \rf{mu}.

We shall construct the counterexample in the setting of Remark \ref{rem-dimension} (i). Set $n=d=1$, so the $\tilde X$ in \rf{XB} is two dimensional. We will use the notation $\tilde \Bx =(\o, \Bx) \in \dbX_2$, where $\o$ and $\Bx$ refer to the paths of $B$ and $X$, respectively.

\begin{example} \rm
\label{eg-main}
Let $A := [0, 1]$, $x_0=0$, $b(t,\tilde\Bx, a) := a$, $\si(t,\tilde\Bx,a) := 1$, for $(t, \tilde\Bx, a)\in [0, T]\times \dbX_2\times A$, namely SDE \rf{Xpath} (or say, the second equation of \rf{XB}) becomes:
\bel{Xeg}
X_t =  \int_0^t \a_s ds + B_t.
\ee
Moreover, $g(\tilde\Bx) := 1_D(\tilde\Bx)$, where, for $\tilde\Bx = (\o, \Bx)\in \dbX_2$,
\bel{D}
\left.\ba{c}
\ds\a^*(t,\tilde\Bx):=0 \vee \Big[\limsup_{h\to0}{(\Bx-\o)_{t}-(\Bx-\o)_{(t-h)^+}\over h}\Big] \wedge 1,\\
\ds D:= \Big\{\tilde \Bx\in \dbX_2: \int_0^T |\a^*(t,\tilde \Bx) - \mu(t, \Bx)|dt=0\Big\}.
\ea\right.
\ee
Then $V^o_0 =0 < 1= V^c_0$.
\end{example}
\proof We first prove $V^c_0 = 1$. Since $g \les 1$, it is clear that $V^c_0\les 1$. Next, by Girsanov theorem, SDE \rf{Xmu} has a unique (in law) weak solution $(\O, \cF, \dbP, B, X)$. Denote $\tilde X := (B, X)$ as usual and consider the closed-loop control $\a_t := \mu(t, X)$, which is obviously $\dbF^{X}$-progressively measurable. Note that $\mu$ takes values in $[0, 1]$. Then $(\O, \cF, \dbF^{\tilde X}, \dbP, B, X, \a)$ is a weak solution to SDE \rf{Xeg} in the sense of Definition \ref{defn-weak} and $\a\in \dbL^0(\dbF^{\tilde X})$. Therefore, $V^c_0 \ges \dbE^\dbP[g(\tilde X)] = \dbP(\tilde X \in D)$. By \rf{Xmu}, it is clear that
$$
\limsup_{h\to0}{(X- B)_t-(X-B)_{(t-h)^+}\over h} = \limsup_{h\to0}{1\over h}\int_{(t-h)^+}^t \mu(s, X_\cd)ds = \mu(t, X_\cd),\q dt\times d\dbP\mbox{-a.s.}
$$
Then
$$
\a^*(t, \tilde X_t) = \mu(t, X_t), ~dt\times d\dbP\mbox{-a.s.},\q\mbox{and thus}\q \tilde X\in D,~\dbP\mbox{-a.s.}
$$
This implies $V^c_0 \ges \dbP(\tilde X\in D) = 1$, and therefore, $V^c_0=1$.

It remains to show that $V^o_0=0$. Let $(\O, \cF, \dbF, \dbP, B, X, \a)$ be an arbitrary weak solution to SDE \rf{Xeg} with open-loop control $\a\in \dbL^0(\dbF^B)$. Note that in this case $X\in \dbL^0(\dbF^B)$ is a strong solution, then $\tilde X\in \dbL^0(\dbF^B;\dbR^2)$.
For the $t_k$ in Lemma \ref{lem-Tsirelson}, introduce:
\bel{Ek}
E_k := \Big\{\int_0^{t_k}  |\a_t - \mu(t, X)|dt=0\Big\},~ k\les -1,\qq E_\infty:= \lim_{k\to -\infty} E_k.
\ee
Note that $E_k \uparrow E_\infty$ as $k\to -\infty$. Clearly $E_k \in \cF^B_{t_k}$, then by the Blumenthal $0$-$1$ law we have $\dbP(E_\infty) = 0$ or $1$. If $\dbP(E_\infty) = 0$, since $\{\tilde X \in D\} = E_0 \subset E_\infty$, then $\dbE^\dbP[g(\tilde X)] = \dbP(\tilde X\in D) = 0$, which is the desired equality we want. So from now on we assume by contradiction that $\dbP(E_\infty) = 1$.

For each $k\les -1$, introduce $\a^k, X^k\in \dbL^0(\dbF^B)$ as follows:
\bel{ak}
\left.\ba{c}
\ds \a^k_t :=  \a_t, ~ t \in [0, t_k);\q \a^k_t := \th\Big({B_{t_i} - B_{t_{i-1}} + \int_{t_{i-1}}^{t_i} \a^k_s ds\over t_i - t_{i-1}}\Big),\q t\in [t_{i}, t_{i+1}),~ i= k,\cds, -1;\\
\ds X^k := \int_0^t \a^k_s ds + B_s.
\ea\right.
\ee
Note that, for $i<k$,
\bel{ai}
\a_t = \mu(t, X)=\th\Big({B_{t_i} - B_{t_{i-1}} + \int_{t_{i-1}}^{t_i} \a_s ds\over t_i - t_{i-1}}\Big),\q  dt\times d\dbP\mbox{-a.s. on}~[t_{i}, t_{i+1})\times E_k.
\ee
Now for $n<k$, since $E_k$ is increasing as $k\to -\infty$, clearly $(\a^n_t, X^n_t) = (\a^k_t, X^k_t) = (\a_t, X_t)$ for $t\les t_n$, and for $i=n,\cds, k-1$, by applying \rf{ak} for $n$ and \rf{ai} for $k$  we see that
\bel{ank}
(\a^n_t, X^n_t) = (\a^k_t, X^k_t),\q \a^k_t = \mu(t, X^k),\q \mbox{on}~ [t_i, t_{i+1})\times E_k.
\ee
Then by applying \rf{ak} for both $n$ and $k$ and recalling \rf{mu} we see that \rf{ank} holds on  $[t_i, t_{i+1})\times E_k$ for $i=k,\cds, -1$ as well. That is, \rf{ank} holds $dt\times d\dbP$-a.s. on  $[0, T]\times E_k$ for all $n\les k$. {\color{black} In particular, this implies the following limits exist:
\bel{hataX}
\hat \a := \lim_{k\to -\infty} \a^k \in \dbL^0(\dbF^B),\q \hat X := \lim_{k\to -\infty} X^k \in \dbL^0(\dbF^B),\q dt\times d\dbP\mbox{-a.s.}
\ee
}
Then, by \rf{ak}  we have
\bel{hatmu}
\hat \a_t = \mu(t, \hat X),\q \hat X_t = \int_0^t \hat \a_sds + B_t,
\ee
  $dt\times d\dbP$-a.s. on $[0, T]\times E_k$ for each $k$, and thus $dt\times d\dbP$-a.s. on $[0, T]\times E_\infty$. By the assumption $\dbP(E_\infty)=1$, we see that \rf{hatmu} holds $dt\times d\dbP$-a.s. on $[0, T]\times \O$. This implies that $\hat X$ is a strong solution of SDE \rf{Xmu}, which is a desired contradiction. So $\dbP(E_\infty)=0$ for all weak solutions with open-loop controls, and therefore $V^o_0 =0$.
  \qed

\begin{remark}
\label{rem-open}
{\rm In this remark we present some related interesting questions we would like to explore in the future research.

(i) The above counterexample relies heavily on the path dependence of the terminal condition $g$, and control only enters in the drift. Is it possible to construct a counterexample such that all the coefficients are state dependent, and/or the control appears in the diffusion as well?

(ii) Regardless whether the open-loop and closed-loop dynamic value functions are equal or not, they might be discontinuous in general. Is it possible to establish the connection between these value functions and the so-called discontinuous viscosity solutions of  the HJB equations?
\qed}
\end{remark}

\section{\color{black}Some further discussions}
\label{sect-analysis}
\setcounter{equation}{0}
In this section we provide some further discussions on the relationship between $V^o_0$ and $V^c_0$ in general setting, without invoking the viscosity solution approach. {\color{black}The arguments are rather standard and the conditions are restrictive in some aspects. Our main point is that these results do not require the continuity of  the coefficients, especially $g$. It will be very interesting to explore more general results when we lose the desired regularity, which we leave to future research.}

\begin{proposition} \sl
\label{prop-Vo<Vc}
Assume $n=d$,  $\si$ takes values in $\dbS^n$ and is positive definite, and $b=b(t, \Bx)$ does not depend on $\a$. Then we have $V^o_0 \les V^c_0$.
\end{proposition}
\proof  Let $(\O, \cF, \dbF, \dbP, B, X, \a)$ be an arbitrary weak solution of \rf{Xpath} such that $\a\in \dbL^0(\dbF^B)$. Note that the quadratic variation process $\langle X\rangle$ is in $\dbL^0(\dbF^X; \dbS^n)$, then so is $\si(t, X, \a_t) = \big({d\over dt} \langle X\rangle_t)^{1\over 2}$, thanks to the assumption that $\si$ is positive definite. Note that
$$
dB_t=\si^{-1}(t,X,\a_t)[dX_t-b(t,X)dt]=\({d\over dt}\lan X\ran_t\)^{-{1\over2}}[dX_t-b(t,X)dt].
$$
Then $B\in \dbL^0(\dbF^X)$ and thus $\a\in \dbL^0(\dbF^B)\subseteq\dbL^0(\dbF^X)$. This implies $\dbE^\dbP[g(X)] \les V^c_0$, hence $V^o_0 \les V^c_0$.
\qed

\ms
Following Krylov \cite{Krylov1}, Gyongy \cite{Gyongy}, and Brunck-Shreve \cite{BS}, we have the following result in the state dependent case.

\begin{proposition} \sl
\label{prop-state}
Assume $b, \si, g$ are state dependent, and for any $(t,x)\in [0, T]\times \dbR^n$, the set $\big\{\big(b(t,x,a), \si\si^\top(t,x,a)\big): a\in A\big\}\subset \dbR^n\times \dbR^{n\times d}$ is convex. Then $V^o_0 \les V^c_0$.
\end{proposition}
\proof  Note that we allow $\si$ to be degenerate. By increasing the dimension of either $B$ or $X$ to $n\vee d$, if necessary, we may assume without loss of generality that $n=d$.

Let $(\O, \cF, \dbF, \dbP, B, X, \a)$ be an arbitrary weak solution of  SDE \rf{Xpath} in the state dependent setting such that $\a\in \dbL^0(\dbF^B)$. By setting $Y:= X-x_0$ and $Z := X$ in \cite[Theorem 3.6]{BS}, we have
\begin{enumerate}[(i)]
\item there exists a measurable function $(\hat b, \hat \si): [0, T]\times \dbR^n \to \dbR^n\times \dbS^{n}$ such that
\bel{hatbsi}
\hat b(t, X_t) = \dbE^\dbP\big[b(t, X_t, \a_t) |X_t\big],\q \hat\si^2(t, X_t) = \dbE^\dbP\big[\si\si^\top(t, X_t, \a_t) |X_t\big];
\ee

\item there exists a probability space $(\hat \O, \hat \cF, \hat\dbP)$, a Brownian motion $\hat B$, and a process $\hat X$ such that
\bel{hatX}
\hat X_t = x_0 + \int_0^t \hat b(s, \hat X_s) ds + \int_0^t \hat \si(s, \hat X_s) d\hat B_s,\q\hat\dbP\mbox{-a.s.}
\ee

\item for any $t$, the $\hat \dbP$-distribution of $\hat X_t$ is equal to the $\dbP$-distribution of $X_t$.
\end{enumerate}

\no Since $\{(b(t,x,a), \si\si^\top(t,x,a)): a\in A\}$ is convex, by \rf{hatbsi} there exists a measurable mapping $\hat \a: [0, T]\times \dbR^n\to A$ such that
\bel{hata}
\hat b(t, X_t) = b(t, X_t, \hat\a(t, X_t)),\q \hat\si^2(t, X_t) = \si\si^\top(t, X_t, \hat \a(t, X_t)).
\ee
Moreover, there exists a mapping $Q: [0, T]\times \dbR^n \to \dbR^{n\times n}$ such that $Q(t,x)$ is an orthogonal matrix and $\hat \si(t,x) = \si(t,x, \hat\a(t,x)) Q(t,x)$. Denote $\tilde B_t := \int_0^t Q(s, \hat X_s) d\hat B_s$, which is still a $\hat \dbP$-Brownian motion. Then \rf{hatX} and \rf{hata} imply
$$
\hat X_t = x_0 + \int_0^t b(s, \hat X_s, \hat\a(s, \hat X_s)) ds + \int_0^t  \si(s, \hat X_s, \hat\a(s, \hat X_s)) d\tilde B_s,\q\hat\dbP\mbox{-a.s.}
$$
This means that $(\hat \O, \hat \cF, \dbF^{\hat B, \hat X}, \hat \dbP)$ and  $(\tilde B, \hat X, \hat \a(\hat X))$ is a weak solution to \rf{Xpath} and $\hat\a(\hat X)\in \dbL^0(\dbF^{\hat X}, \hat\dbP)$, and thus $\dbE^{\hat\dbP}[g(\hat X_T)] \les V^c_0$. Finally, by (iii) we have $\dbE^\dbP[g(X_T)] = \dbE^{\hat\dbP}[g(\hat X_T)] \les V^c_0$.
\qed

\ms
We remark that, in (iii) above, only the marginal distributions are equal. In general the $\hat \dbP$ (joint) distribution of the process $\hat X_{[0, T]}$ does not coincide with the $\dbP$ distribution of $X_{[0, T]}$, so we are not able to extend these arguments to the path dependent case.

\begin{proposition}
\label{prop-Vc<Vo} \sl
Under the following two conditions we have $V^c_0 \les V^o_0$:

{\rm(i)} $\si=\si(t, \Bx)$ does not depend on $\a$ and is uniformly Lipschitz continuous in $\Bx$;

{\rm(ii)} $b = \si \l$ where the function $\l: [0, T]\times \dbX_n \times A\to \dbR^d$ is bounded, continuous in $a$, and uniformly Lipschitz continuous in $\Bx$.
\end{proposition}
\proof  Let $(\O, \cF, \dbF, \dbP, B, X, \a)$ be an arbitrary weak solution of \rf{Xpath} such that $\a\in \dbL^0(\dbF^X)$.  It suffices to show that $\dbE^\dbP[g(X_\cd)] \les V^o_0$. For this purpose, we denote
$$
B^\a_t:= B_t + \int_0^t \l(s, X_\cd, \a_s)ds,\q {d\dbP^\a\over d\dbP} := M^\a_T := e^{-\int_0^T \l(s, X_\cd, \a_s) dB_s -{1\over 2} \int_0^T |\l(s, X_\cd, \a_s)|^2 ds}.
$$
By Girsanov Theorem, we know that $\dbP^\a\sim\dbP$ (meaning that they are equivalent) and $B^\a$ is a $\dbP^\a$-Brownian motion. Note that
\bel{Xa}
X_t = x_0 + \int_0^t \si(s, X) dB^\a_s.
\ee
Fix a probability space $(\O^0, \cF^0, \dbP^0)$ and a Brownian motion $B^0$ on it. Under (i) the SDE
\bel{X0}
X^0_t = x_0 + \int_0^t \si(s, X^0_\cd) dB^0_s, \q\dbP^0\mbox{-a.s.}
\ee
 has a unique strong solution $X^0$.  Now compare \rf{Xa} and \rf{X0} we see that the $\dbP^\a$-distribution of $(B^\a, X)$ is equal to the $\dbP^0$-distribution of $(B^0, X^0)$. Since $\a\in \dbL^0(\dbF^X, \dbP)$, we may write it as $\a(t, X_\cd)$. Then
\bel{Girsanov}
\left.\ba{c}
\ds \dbE^\dbP[g(X_\cd)] = \dbE^{\dbP^\a}\big[(M^\a_T)^{-1} g(X_\cd)\big] = \dbE^{\dbP^0}\big[N^\a_T g(X^0_\cd)\big],\\
\ns\ds\mbox{where}\q N^\a_T:= \exp\Big(\int_0^T \l(s, X^0_\cd, \a(s, X^0)) dB^0_s -{1\over 2} \int_0^T |\l(s, X^0_\cd, \a(s, X^0))|^2 ds\Big).
\ea\right.
\ee

For fixed $\dbP^0$,  there exist piecewise constant processes $\ds\a^n(t, X^0) = \sum_{i=0}^{n-1} \a^n(t_i, X^0) 1_{[t_i, t_{i+1})}(t)$ such that $\lim_{n\to\infty} \dbE^{\dbP^0}\big[\int_0^T |\a^n(t, X^0) - \a(t, X^0)|^2dt \big]=0$.   Then by (ii) one can easily show that
\bel{Nan}
\lim_{n\to\infty} \dbE^{\dbP^0}\big[|N^{\a^n}_T - N^\a_T|^2\big]=0,\q\mbox{and hence}\q \lim_{n\to\infty}\dbE^{\dbP^0}\big[N^{\a^n}_T g(X^0_\cd)\big]=\dbE^{\dbP^0}\big[N^\a_T g(X^0_\cd)\big].
\ee
For each $n$, by the Girsanov theorem we have $\dbE^{\dbP^0}\big[N^{\a^n}_T g(X^0_\cd)\big] = \dbE^{\dbP^n}\big[g(X^0)\big]$, where $\dbP^n\sim \dbP^0$ is a probability measure, $B^n_t := B^0_t - \int_0^t \l(s, X^0, \a^n(s, X^0))ds$ is an $\dbP^n$-Brownian motion, and
$$
X^0_t = x_0 + \int_0^t b(s, X^0, \a^n(s, X^0)) ds + \int_0^t \si(s, X^0)dB^n_s.
$$
Since $\a^n$ is piecewise constant, and $b, \si$ are uniformly Lipschitz continuous in $\Bx$, by induction on $i$ one can easily show that $\dbF^{X^0} \subset \dbF^{B^n}$, where the augmentation is under $\dbP^0$ and equivalently under $\dbP^n$. Then $\a^n(t, X^0)\in \dbL^0(\dbF^{B^n}, \dbP^0)$, namely is an open-loop control. This implies that $\dbE^{\dbP^0}\big[N^{\a^n}_T g(X^0_\cd)\big] = \dbE^{\dbP^n}\big[g(X^0)\big] \les V^o_0$. Then by \rf{Girsanov} and  \rf{Nan} we have $\dbE^\dbP[g(X_\cd)] = \dbE^{\dbP^0}\big[N^\a_T g(X^0_\cd)\big] \les V^o_0$, and therefore $V^c_0\les V^o_0$.
\qed

\ms

Note again that in the above proof,  we may allow $g$ to be discontinuous. Combine Propositions \ref{prop-state} and \ref{prop-Vc<Vo}, we immediately have the following.

\begin{corollary} \sl
\label{cor-Vc=Vo}
Assume  $b, \si, g$ are state dependent, and

{\rm(i)} $\si=\si(t, x)$ does not depend on $\a$ and is uniformly Lipschitz continuous in $x$;

{\rm(ii)} $b = \si \l$ where the function $\l: [0, T]\times \dbR^n \times A\to \dbR^d$ is bounded, continuous in $a$, uniformly Lipschitz continuous in $x$, and the set $\{\l(t,x,a): a\in A\}\subset \dbR^d$ is convex.

Then $V^o_0 = V^c_0$.
\end{corollary}

\begin{remark}
\label{rem-Vc=Vo}
{\rm
Under the conditions in Corollary \ref{cor-Vc=Vo}, obviously the dynamic value functions for the control problem on $[t, T]$ with initial value $x$ are also equal:  $v^o(t,x) = v^c(t,x) =: v(t,x)$. However, we emphasize here that $g$ can be discontinuous and $\si$ can be degenerate, then $v$ might be discontinuous. One trivial example is: $b=0, \si=0$, then $v(t,x) = g(x)$ for all $t$, which will be  discontinuous if $g$ is so.
\qed}
\end{remark}

 \end{document}